\documentclass[12pt]{article}
\usepackage{amsmath,amsfonts,amssymb,amsthm}

\def\noi{\noindent}

\def\pf{\noi{\bf Proof.\ \,}}
\def\eop{{$\square$}}

\def\labtt#1{\label {#1}}

\def\a{\alpha}
\def\b{\beta}
\def\g{\gamma}

\def\s{\sigma}

\def\CC{{\mathbb C}}

\def\RR{{\mathbb R}}

\def\ZZ{{\mathbb Z}}

\def\la{\langle}
\def\ra{\rangle}
\def\<{\langle}
\def\>{\rangle}

\def\bs{\it}            

\def\dim{{\bs dim}}

\def\vac{\hbox{\bf 1}} 

\def\veehp{V_{EE_8}^+}

\def\mcb{{\mathcal B}}

\def\ordermonster{2^{46}3^{20}5^97^611^213^317{\cdot}19{\cdot}23{\cdot}39{\cdot}31{\cdot}41{\cdot}47{\cdot}59{\cdot}71}

\def\cvcch{cvcc\frac 12}

\begin{document}

\newtheorem{thm}{Theorem}[section]
\newtheorem{prop}[thm]{Proposition}
\newtheorem{lem}[thm]{Lemma}
\newtheorem{rem}[thm]{Remark}
\newtheorem{coro}[thm]{Corollary}
\newtheorem{conj}[thm]{Conjecture}
\newtheorem{de}[thm]{Definition}
\newtheorem{hyp}[thm]{Hypothesis}

\newtheorem{nota}[thm]{Notation}
\newtheorem{ex}[thm]{Example}
\newtheorem{proc}[thm]{Procedure}

\def\refp#1{(\ref {#1})}
\def\refpp#1{(\ref {#1})}

\begin{center}
{\Large  \bf
A new existence proof of the Monster by VOA theory }

\bigskip

\centerline{March 7, 2011}

\vspace{10mm}
Robert L.~Griess Jr.
\\[0pt]
Department of Mathematics\\[0pt] University of Michigan\\[0pt]
Ann Arbor, MI 48109  USA  \\[0pt]
{\tt rlg@umich.edu}\\[0pt]
\vskip 5mm

Ching Hung Lam
\\[0pt]
Institute of Mathematics \\[0pt]
Academia Sinica\\[0pt]
Taipei 10617, Taiwan\\[0pt]
{\tt chlam@math.sinica.edu.tw}\\[0pt]
\vskip 1cm

\end{center}

\begin{abstract}
We use uniqueness of a VOA (vertex operator algebra) extension of  $(V_{EE_8}^+)^3$
to a Moonshine type VOA to give a new existence proof of a finite simple group
of Monster type. The proof is relatively direct. Our methods depend on VOA
representation theory and are free of many special calculations which
traditionally occur in theory of the Monster.
\end{abstract}

\newpage

\tableofcontents

\section{Introduction}

We define a finite group $G$ to be of {\it Monster type} if it has an   involution
$z$ whose centralizer $C_G(z)$  has
the form $2^{1+24}Co_1$, is 2-constrained (i.e., satisfies $\la z \ra = C_G(O_2(C_G(z)))$ and $z$ is conjugate to an element in
$C_G(z)\setminus \{ z\}$.  A short argument proves that such a $G$ must
be simple (e.g., see \cite{grfg,tits84}).
The present article gives a new and relatively direct
existence proof of a group of Monster type.
In fact, a group of Monster type is unique up to isomorphism \cite{gms},
so the group constructed in this article can be called ``the'' Monster, the group constructed in \cite{grfg}.
To avoid specialized finite group theory in this article,
we work with a group of Monster type and refer to \cite{gms} for uniqueness.

Our basic strategy is described briefly in the next paragraph.
It was inspired by the article of
Miyamoto \cite{M}, which showed how to make
effective use of {\it simple current modules and extensions}.
Later in this introduction, we sketch these important concepts.
In a sense, our existence proof is quite short.
The hard group theory and case-by-case analysis of earlier proofs have  essentially been eliminated.
We use the abbreviation VOA  for  {\it vertex operator algebra} \cite{flm}.
Most of this article is dedicated to explaining how existing VOA theory applies.

In  \cite{sh}, Shimakura takes
$(\veehp )^3$ and builds a candidate $V$  for the Moonshine VOA (a short account is in Section
\ref{sec:2.1} of the present article).  His construction furnishes a large subgroup  of
$Aut(V)$.  From this subgroup, we take a certain involution and analyze $V^+, V^-$,
its fixed point VOA and its negated space on $V$, respectively.  We can recognize $V^+$ as a Leech
lattice type VOA.  The group $Aut(V^+)$  and its extension (by projective representations) to irreducibles of the fixed point VOA
are understood.  One of these irreducibles is $V^-$.
We thereby get a new subgroup of $Aut(V)$ which has the shape $2^{1+24}Co_1$ and is moreover isomorphic to the
centralizer of a 2-central involution in the Monster.  These two subgroups of $Aut(V)$ generate the larger group $Aut(V)$, which we then prove is a finite group  of Monster type.

We  refer to the overVOA   of
$(V_{EE_8}^+)^3$ constructed in \refpp{V}  as a VOA of {\it Moonshine type},
meaning a holomorphic VOA $V=\oplus_{n=0}^{\infty} V_n$ of central charge 24, so that $V_0$ is 1-dimensional, $V_1=0$ and the Monster acts as
a group of automorphisms with faithful action on $V_2$.
We mention that such a VOA is
isomorphic to the standard Moonshine VOA constructed in \cite{flm}, by
\cite{dgl,ly}.
For the purpose of this article, it is not necessary to quote such characterizations.

The theory of simple current modules originated in
the papers \cite{FG} and \cite{SchY}. In \cite{dlm,li}, certain simple current modules of a VOA are constructed using weight one semi-primary elements and extensions of a VOA by its simple current modules are also studied. The notion of simple current extension turns out to be a very powerful tool for constructing new VOAs from a known one  \cite{dgh,lly,ly1,li,M96,M}.
Let $V$ be a VOA and let $\mathcal{M}=\{M^i\mid i\in I\}$ be a finite set of irreducible modules of $V$ with integral weights. If $V\in \mathcal{M}$ and $\oplus_{M\in \mathcal{M}} M$ is closed under the fusion rules, it is
possible that $\oplus_{M\in \mathcal{M}} M$  carries the structure of a VOA for which $V$ is a subVOA.
In general,  it is extremely difficult to determine if $\oplus_{M\in \mathcal{M}} M$ has  a VOA structure.   There may be no such VOA structures, or there could be many.
When the simple current property holds \refpp{SCE}, there is  a VOA
structure on $\oplus_{M\in \mathcal{M}} M$ extending the given action of $V$ and
the VOA structure is unique if the underlining field is algebraically closed
\cite[Proposition 5.3]{dm} (see also \cite{dmz,M96}).
This ``rigidity'' of simple current extensions is useful in structure
analysis and leads to certain transitivity results, which reduce the need for  calculations.

As described in \cite{grfg}, existence of the Monster implies existence of several other sporadic groups which had originally been constructed with special
methods, including computer work.   We hope that the present
article may suggest  useful viewpoints for other sporadic groups.

\subsection{About existence proofs of the Monster}

The first existence proof of the Monster was made in 1980, and published in \cite{grfg}; see also \cite{grfgpnas}.   A group $C\cong 2^{1+24}Co_1$ and a representation of degree 196883 was described.  The hard part was to choose a $C$-invariant algebra structure, give an automorphism $\s$ of it which did not come from $C$, then identify the group $\la C, \s \ra$ by proving finiteness and proving that  $C$ is an involution centralizer in it.

During the decade  that followed \cite{grfg}, there were analyses, improvements
and alternate viewpoints by Tits  \cite{tits82a,tits82b,tits83,tits84}  and
Conway \cite{conway}. In the mid-80s, the theory of vertex algebras was developed. The
Frenkel-Lepowsky-Meurman
text \cite{flm} established the important construction of a Moonshine VOA and
became a basic reference for VOA theory.  The construction of the Monster done
in \cite{flm} followed the lines of \cite{grfg}, but in a broader VOA setting. The articles
\cite{dgm,dgm2} constructed a VOA and gave a physics field theory
interpretation to aspects of \cite{flm,grfg}.

In 2004,  Miyamoto \cite{M},
made significant use of simple current extensions to give a new construction of a Moonshine VOA and of the Monster acting as automorphisms.
An existence proof of the Monster was recently announced in \cite{ivanov},
which uses theories of finite geometries and group amalgams.

Uniqueness was first proved in \cite{gms}.  A different uniqueness proof is indicated in \cite{ivanov}.

\subsection{Acknowledgments}

The first author thanks the United States National Science Foundation (NSF DMS-0600854)
and National Security Agency (NSA H98230-10-1-0201)
for financial support
and the Academia Sinica for hospitality during a visit August, 2010.

The second author thanks  National Science Council (NSC 97-2115-M-006-015-MY3)
and National Center for Theoretical Sciences, Taiwan for financial support.

\bigskip

\newpage
\subsection{Table of Notations }
{\small
\begin{tabular}{|c|c|c|}
   \hline
\bf{Notation}& \bf{Explanation} & \bf{Examples in text}  \cr
   \hline\hline
   $A.B, A{:}B,$ & group extension of normal subgroup $A$ by quotient& \refpp{normalizertrio}  \cr
  $A{\cdot}B$  & $B$, split extension, nonsplit extension, respectively &\refpp{c=cz}, \cr \hline
$BW_{16}$ & the Barnes-Wall lattice of rank $16$ & Lemma \ref{BW}\cr \hline
$C$ & the centralizer of an involution in $Aut(V)$ & \refpp{bigc}  \cr\hline
$\cvcch$   & simple conformal vector  & Sec. \ref{sec:3.1}\cr
& of central charge $\frac{1}2$ & \cr \hline
$Co_1$ & the first Conway group $O(\Lambda)/\{\pm 1\}$ & Page 2 \cr \hline
$EE_8$ & lattice isometric to $\sqrt 2$ times & Page 2 \cr
&   the famous  $E_8$ lattice & Sec \ref{sec:2.1}\cr
\hline
$\mathbb{F}\{L\}$ & a twisted group algebra of a lattice $L$ over a field $\mathbb{F}$ & Coro. \ref{VLeech}\cr \hline
$g\circ M$ & $g$-conjugate module of a $V$-module $M$ & Notation \ref{gconj} \cr
\hline
$H_\Lambda$ & the subgroup of  $Aut(V_\Lambda^+)$   generated  & Notation \ref{HL}\cr
&by Miyamoto involutions of $AA_1$-type       &          \cr \hline
  & the subgroup of  $Aut(V_\Lambda^+)$   generated  & Notation \ref{HL}\cr
${\tilde H}_\Lambda$ &by Miyamoto involutions of $AA_1$-type       &          \cr \hline
$\Lambda$ & Leech lattice, the unique even unimodular  & Sec. \ref{C(z)}\cr
& lattice of rank $24$ with no roots& Coro. \ref{VLeech}\cr \hline
$M(1)$ & the unique irreducible
   $\hat{\mathfrak{h}}$-module such that &  Coro. \ref{VLeech}\cr
    & $\alpha\otimes t^n 1=0$ for all $\alpha\in {\mathfrak h}$ and $n>0$ and $K=1$,& \cr & where $\mathfrak{h}= \CC\otimes_\ZZ L$ and $\hat{\mathfrak{h}}= \oplus_{n\in \ZZ} (\mathfrak{h}\otimes t^n)\oplus \CC K$  & \cr \hline
$M\times_V N$ & the fusion product of $V$-modules $M,N$ & Notation \ref{fusionP}\cr \hline
$O_p(G)$ & maximal normal $p$-subgroup of $G$ & Page 2 \cr \hline
$O_{p'}(G)$ & maximal normal subgroup of $G$ of order prime to $p$  & Page 2 \cr\hline
$R(U)$ & the set of all inequivalent irreducible modules of $U$ & Sec. \ref{sec:2.1}\cr\hline
$Stab_G(X)$ & subgroup of the group $G$ which stabilizes the set $X$& \refpp{framestab} \cr \hline
 $t(e)$ & Miyamoto involution associated  to a $\cvcch$ $e$ & Notation \ref{MI} \cr \hline
 $U$ & the VOA $\veehp$ & \refpp{u} \cr \hline
 $2^{1+2n}$ & extra-special 2-group of order $2^{1+2n}$ & Page 2\cr \hline
 $2^{1+24}Co_1$ & an extension of $Co_1$ by $2^{1+24}$ & Page 2, Sec. \ref{C(z)}\cr \hline
  $V$ & a VOA which is a simple current extension of $U$ & \refpp{VV}  \cr \hline
 $V_L$ & lattice VOA for positive definite even lattice $L$ & Lem. \ref{V124}, Coro. \ref{VLeech} \cr \hline
  $V_{L, \RR} $ & lattice VOA over $\RR$  for positive definite  & Lem. \ref{V124}, Coro. \ref{VLeech} \cr
  &even lattice $L$ &\cr \hline
  &  a real form of the lattice VOA  $V_L$ whose  &Prop.  \ref{VLpd} \cr
  $\tilde{V}_{L,\RR}$  &  invariant form is positive definite& \cr\hline
 $V_L^+$&  the fixed point subVOA of $V_L$ & Sec. \ref{sec:2.1} \cr
        &   by a lift of the $(-1)$-isometry of $L$ & Coro. \ref{VLeech}\cr \hline
 \end{tabular}
}

\section{Simple current extensions}
In this section, we shall recall the notion of simple current extensions and their basic properties \cite{dlm,sh}.

Let $V$ be a VOA and let $M_1$, $M_2$, $M_3$ be $V$-modules.
Let $I_V\binom{M_3}{M_1\ M_2}$ be the space
of all $V$-intertwining operators of type
$\binom{M_3}{M_1\ M_2}$. We refers to \cite[Chapter 5]{FHL} for the definition of intertwining operators.

The {\it fusion coefficient (or fusion rule)} is defined to be the integer
\[
N_{M_1\ M_2}^{M_3} =\dim \left( I_V\binom{M_3}{M_1\quad M_2}\right).
\]

\begin{nota}\labtt{fusionP}
Suppose $V$ is rational. Then $V$ has only finitely
many inequivalent irreducible modules and all $V$-modules are direct sums of finitely
many irreducibles. Moreover, all fusion coefficients are finite \cite{FZ}.
We define a formal product of two isomorphism types of irreducible modules $M_1$ and $M_2$:
\[
M_1\times_V M_2=\sum_{M_3}  N_{M_1\ M_2}^{M_3} M_3,
\]
where $M_3$ runs through the finitely many isomorphism types of irreducible modules of $V$. This product is often called
the fusion product.
\end{nota}

\begin{de}\cite{dlm}\labtt{SC}
Let $V$ be a rational VOA. A $V$-module $M$ is called a \textit{simple current module} if
the fusion product $M\times_V N$ is again irreducible for any irreducible $V$-module $N$, i.e.,
\[
\sum_{W \text{ irred}}  N_{M\ N}^{W} =1.
\]
By definition, it is clear that $M\cong M\times_V V$ is irreducible if $M$ is a simple current
module.
\end{de}

\begin{de}\labtt{fullsubvoa}
A {\it full} subVOA is a subVOA which contains the principal Virasoro element of the larger VOA.
\end{de}

Now let $V^0$ be a simple rational $C_2$-cofinite VOA of CFT type and
let $D$ be a finite abelian group. Let $\{
V^\alpha \mid \alpha \in D\}$ be a set of inequivalent irreducible $V^0$-modules
indexed by $D$. Assume that the weights of $V^\a, \a\in D,$ are integral and $$V^\a \times_{V^0} V^\b = V^{\a+\b}$$
for all $\a, \b \in D$.

\begin{de}[\cite{dlm}] \labtt{SCE}
A VOA $V$  is called a {\it ($D$-graded)
simple current  extension} of $V^0$ if

1.  $V^0$ is a full subVOA of $V$ ;

2.  $V=\oplus_{\alpha\in D} V^\alpha$;

3.  All $V^\alpha$, $\alpha\in D$, are simple current
$V^0$-modules.

\end{de}

\medskip

Next we shall recall the notion of $g$-conjugate modules \cite{dm0,sh}.

\begin{de}\labtt{gconj}
Let $(M, Y_M)$ be a $V$-module and $g\in Aut(V)$ an automorphism.  The $g$-conjugate module of $M$
is defined to be the $V$-module $(g\circ M, Y_{g\circ M})$, where $g\circ M=M$ as a vector space and
$Y_{g\circ M}(v,z)= Y_M(g^{-1} v, z) $ for all $v\in V$.
\end{de}

\begin{rem}\labtt{nu}
 By definition, there exists a linear isomorphism $\nu: g\circ M \to M$ such that
$\nu Y_{g\circ M}(v,z)= Y_M(g^{-1} v, z)\nu $  for any $v\in V$. In fact,  one can  assume
$\nu=id_M$ by identifying $g \circ M$ with $M$ as  vector spaces.
\end{rem}

The following lemma can be shown easily by definition.

\begin{lem}[cf. \cite{sh}]
Let $V$ be a VOA and $g\in Aut(V)$.
\begin{enumerate}
\item If $W$ is an irreducible $V$-module, then $g\circ W$ is also irreducible.

\item  If $M\cong W$ as a $V$-module, then $g\circ M\cong g\circ W$.

\item  Let $M_1,M_2, M_3$ be $V$-modules. For each $x\in \{1,2,3\}$,  let $\nu_x: g\circ M_x \to
M_x$ be a linear isomorphism such that $\nu_x Y_{g\circ M_x}(v,z) =
Y_{M_x}(g^{-1}v,z) \nu_x$ for any $v\in V$. Then for any nonzero intertwining operator $\mathcal{Y}$
of type $\binom{M_3}{M_1\ M_2}$, the map $ \nu_3 \mathcal{Y}( \nu_1 u,z)\nu_2$ is a nonzero
intertwining operator of type $\binom{g\circ M_3}{g\circ M_1 \ \ g\circ M_2}$. In particular, the
fusion rules are preserved by $g$-conjugation.
\end{enumerate}
\end{lem}

\medskip

The following theorem follows easily by the fusion rules $V^\a \times V^\b =V^{\a+\b}$.

\begin{thm}\labtt{charactersofd}
Let $V^0$ be a rational $C_2$-cofinite VOA of CFT type and let $V=\oplus_{\a\in D} V^\a$  be a
(D-graded) simple current extension of $V^0$. Let $D^*$ be the group of all irreducible characters
of $D$. Then for any $\chi\in D^*$, the linear map
\[
\tau_\chi(v) = \chi(\a) v \qquad \text{ for  any }  v\in V^\a, \a\in D,
\]
defines an automorphism of $V$. In particular, $\{\tau_\chi\mid \chi\in D^*\}\cong D^*$ is an abelian subgroup of $Aut(V)$.
\end{thm}

\begin{nota}\labtt{D*}
By abuse of notation, we often denote the group $\{\tau_\chi\mid \chi\in D^*\}< Aut(V)$ defined in \refpp{charactersofd} by $D^*$.
\end{nota}

\begin{nota}\labtt{Wclass}
Let $W$ be an irreducible $V$-module. We shall use $[W]$ to denote the
isomorphism class containing $W$ .
\end{nota}

The next theorem gives a criterion for lifting an automorphism of $V^0$ to $V$
and can be proved using the general arguments for simple current extensions
\cite{SY,sh1}.

\begin{thm}[cf. \cite{sh1}]\labtt{liftofsce}
Let $V^0$ be a rational $C_2$-cofinite VOA of CFT type and $V=\oplus_{\a\in D}
V^\a$ a (D-graded) simple current extension of $V^0$. Let $g\in Aut(V^0)$. Then
there exists an automorphism $\tilde{g}\in Aut(V)$ such that
$\tilde{g}|_{V^0}=g$ if and only if $\{[g\circ V^\a]\mid \a\in D\} = \{[V^\a]\mid
\a\in D\}$.
\end{thm}

\medskip

\begin{nota}
  Let $\mathcal{X}$ be a set of isomorphism types of irreducible modules of a VOA $V$. We denote
\[
 g\circ \mathcal{X} =\{[ g\circ W]\mid [W]\in \mathcal{X}\}
\]
 for any $g\in Aut(V)$.
\end{nota}

The next theorem follows easily from Theorem \ref{liftofsce} and a proof can be found in \cite{sh}.

\begin{thm}[Corollary 2.2 of \cite{sh2}]\labtt{ND*}
Let $V=\oplus_{\a\in D} V^\a$ a (D-graded) simple current extension of $V^0$. Denote
\[
N_D=\left \{g\in Aut(V^0)\mid \{[g\circ V^\a] \mid \a\in D \}=\{[V^\a]\mid \a \in D\} \right\}.
\]
Then there exists an exact sequence
\[
1\longrightarrow  D^* \longrightarrow N_{Aut(V^0)}(D^*) \overset{\eta} \longrightarrow  N_D \longrightarrow 1,
\]
where $\eta$ is the restriction map to $V^0$ and $D^*$ is identified with the group
$\{\tau_\chi\mid \chi\in D^*\}$.
\end{thm}

\subsection{Simple current extension of $(V_{EE_8}^+)^3$}\labtt{sec:2.1}

In this section, we shall recall a construction of the Moonshine VOA by Shimakura \cite{sh}. First, we sketch an outline of the Shimakura construction.

It is known \cite{sh1} (see also \cite{ly}) that all irreducible modules of $V_{EE_8}^+$ are simple
current modules and the fusion group $R(\veehp )$ of $\veehp$ has a quadratic spaces over $\ZZ_2$.
Moreover, the automorphism group $Aut(\veehp)\cong O^+(10,2)$ acts faithfully on $R(\veehp)$. By the
general theory of simple current extensions, one knows that holomorphic (simple current) extensions
of $(\veehp)^3$ correspond to maximal totally singular subspaces of $R(\veehp)^3$ up to the action
of $Aut( (\veehp)^3)\cong O^+(10,2)\wr Sym_3$. The main idea of Shimakura is to construct a maximal
totally singular subspace  of $R(\veehp)^3$ with large minimal weight and compute certain
automorphism subgroups using the standard theory of quadratic spaces.

Next we shall review some basic properties of the lattice type VOA $V_{EE_8}^+$
\cite{g156,sh1}.

\begin{nota}\labtt{r(u)}
Let $R(U)$ be the set of all inequivalent irreducible  modules of a
VOA $U$. If $U= V_{EE_8}^+$, then $R(V_{EE_8}^+)$ forms a 10-dimensional
quadratic space over $\ZZ_2$ with respect to the fusion rules and the quadratic
form
\begin{equation}
 q([M])=
\begin{cases}
0 & \text{ if the weights of $M$ are in } \ZZ,\\
1 & \text{ if the weights of $M$ are in } \frac{1}2+\ZZ.
\end{cases}
\end{equation}
We shall denote the corresponding bilinear form by $\la \, , \, \ra$.
\end{nota}

Recall that $Aut(V_{EE_8}^+)\cong O^+(10,2)$ and $Aut(V_{EE_8}^+)$ acts as on $R(V_{EE_8}^+)$ as a
group of isometries \cite{g156,sh1}.

\begin{nota}\labtt{u}
From now on, we use $U$ to denote the VOA $V_{EE_8}^+$. and we use $U^n$ to denote the tensor product of $n$ copies of $U$.
\end{nota}

The proof of the following
proposition can be found in \cite{sh1}.
\begin{prop}\labtt{weight}
The group $Aut(U)\cong O^+(10,2)$ acts transitively on non-zero singular
elements and non-singular elements of $R(U)$, respectively.

1. If $[W]$ is a non-zero singular element in $R(U)$, then the minimal weight of
the irreducible module $W$ is $1$ and $\dim (W_1)=8$.

2. If $[W]$ is  a non-singular element, then the minimal weight of $W$ is $1/2$
and $\dim (W_{\frac{1}2}) =1$.
\end{prop}

\begin{nota}\labtt{mathcals}
Since $R(U^3)\cong R(U)^3$, we shall view $R(U^3)$ as  a direct sun of quadratic
spaces.  See \refpp{r(u)} for notations.
The quadratic form and the associated bilinear form are given by
\[
\begin{split}
q(a,b,c)&= q(a)+q(b)+q(c)\quad  \text{ and }\\
\la (a,b,c), (a',b',c')\ra &=\la a,a'\ra+\la b,b'\ra+\la c,c'\ra
\end{split}
\]
for  $ (a,b,c)\in R(U)^3.$

 Following the analysis of \cite{sh}, let $\Phi$ and $\Psi$
be maximal totally singular subspaces of $R(U)$
such that $\Phi\cap
\Psi =0$. Then the space
\begin{equation}\label{S}
\mathcal{S}:=span_{\ZZ_2}\{ (a,a,0), (0,a,a), (b,b,b)\mid a\in \Phi, b\in \Psi\}
\end{equation}
is a maximal totally singular subspace of $R(U)^3\cong R( U^3)$.
\end{nota}

\begin{de}\labtt{minW}
Let $W$ be an irreducible module. We define the minimal weight of $[W]$ to be
the minimal weight of $W$.
\end{de}

\begin{lem}\labtt{V1}
Let $[W]\in \mathcal{S}$.  Then the minimal weight of $[W]$ is $\geq 2$. If the
minimal weight of $[W]$ is $2$, then up to a permutation of the 3 coordinates,
$[W]$ has the form:

1. $(a, a, 0)$, where $a\in \Phi\setminus \{0\}$ or

2. $(a+b+c, a+c, b+c)$, where $a,b\in \Phi\setminus \{0\}$, $c\in \Psi$, $a+c,b+c$ are non-singular
and $a+b+c$ is non-zero singular.
\end{lem}

\pf  First we note that a general element in $\mathcal{S}$ has the form
\[
 (a+b+c, a+c, b+c),
\]
where $a,b\in \Phi$ and $c\in \Psi$.

Since $\Phi\cap \Psi=0$, $a+c=0$ implies that $a=c=0$ if $a\in \Phi$, $c\in \Psi$.
Thus, up to permutations of the 3 coordinates,  a nonzero element in $\mathcal{S}$ has the form

1. $(a,a,0)$, where $a\in \Phi\setminus \{0\}$ or

2. $(a+b+c, a+c, b+c)$, where $a,b\in \Phi\setminus \{0\}$ and $c\in \Psi$.

In Case 1, the minimal weight is 2. In Case 2, the minimal weight is $> \frac{3}2$ and thus $\geq
2$.

If $(a+b+c, a+c, b+c)$, $a,b,c$ nonzero, has the minimal weight $2$, then two of the $a+b+c,a+c,
b+c$ have minimal
weights $1/2$ and the remaining one has minimal weight $1$. Therefore, two are non-singular and one
is a nonzero singular.
\eop

\medskip

By Lemma \ref{V1}, we have the following theorem.

\begin{thm}[Theorem 4.10 of \cite{sh}]\labtt{V}
Let $V:=\mathcal{V}(\mathcal{S}) = \oplus_{[W]\in \mathcal{S}} W $. Then $V$ is a holomorphic
framed VOA of central charge $24$ and $V_1=0$.
\end{thm}

\begin{nota}\labtt{VV}  For the rest of this article, $V$ shall denote the VOA
$V:=\mathcal{V}(\mathcal{S}) = \oplus_{[W]\in \mathcal{S}} W $
of \refpp{V}.
\end{nota}

\begin{rem}
It was also shown in \cite{sh} that the singular space $\mathcal{S}$ defined in
\eqref{S} is the unique (up to $Aut(U^3)$) maximal totally singular subspace of
$R(U^3)$ such that $\mathcal{V}(\mathcal{S})= \oplus_{[W]\in \mathcal{S}} W $
has trivial weight one subspace.
\end{rem}

\begin{thm}[\cite{sh}]
$Aut(V)$ acts transitively on the set of all subVOAs of $V$ which are isomorphic to $U^3$.
\end{thm}

Now let $H$ be the stabilizer of $\Phi$ in $Aut(U)$ and $K$ the stabilizer of $\Psi$
in $H$. Then $H=2^{10}: L_5(2)$ and $K=L_5(2)$ (cf. \cite[Proposition 2.5]{sh}).
Set $H_{(i,j)} =\{ (g_1, g_2, g_3) \mid g_i\in H,  g_i=g_j\text{ and } g_k=1, \text{
for } k\neq i,j\}$, where $i\neq j$ and $i,j \in  \{1,2,3\}$ and
$K_{(1,2,3)}=\{(g,g,g)\mid g\in K\}$.

\begin{lem}[Proposition 2.5 of \cite{sh}]
Let $S$ be defined as in \eqref{S}. Then the stabilizer $N_\mathcal{S}$ of
$\mathcal{S}$ in $Aut(U^3)$ is generated by $O_2(H_{(1,2)}), O_2(H_{(1,3)})$,
$K_{(1,2,3)}$ and $Sym_3$ and it has the shape $2^{20}: ( L_2(5)\times Sym_3)$.
\end{lem}

As a consequence, we have
\begin{prop}[Corollary 4.18 of \cite{sh}]\labtt{normalizertrio}
Let $\mathcal{S}$ be defined as in \eqref{S}. We identify $\mathcal{S}^*$ with the subgroup $\{\tau_\chi\mid \chi \in
\mathcal{S}^*\}$ of $Aut(V)$ which was defined in \refpp{charactersofd}.   Then
$$N_{Aut(V)}(\mathcal{S}^*)=Stab_{Aut(V)} (U\otimes U\otimes U)\cong
2^{15}(2^{20}:(L_5(2)\times Sym_3)).$$
\end{prop}

Set $\mathcal{S}_\Phi =span\{(a,a,0), (0,a,a)\mid a\in \Phi\}$. Then
$$V_\Phi=\mathcal{V}(\mathcal{S}_\Phi) =\oplus_{[W]\in \mathcal{S}_\Phi}W$$
is a subVOA of $V$.

For any coset $x+\mathcal{S}_\Phi\in \mathcal{S}/\mathcal{S}_\Phi$, the subspace
$\mathcal{V}(x+\mathcal{S}_\Phi)=\oplus_{[W]\in x+\mathcal{S}_\Phi}W$ is an irreducible module of
$V_\Phi$. In fact, $\mathcal{V}(x+\mathcal{S}_\Phi)$ is a simple current module of
$V_\Phi$.

\begin{prop}
Let $V= \oplus_{x+\mathcal{S}_\Phi\in \mathcal{S}/\mathcal{S}_\Phi}  \mathcal{V}(x+\mathcal{S}_\Phi)$ be a simple current extension of $V_\Phi$. Then $N_{Aut(V)}((\mathcal{S}/\mathcal{S}_\Phi)^*)$ stabilizes the subVOA $U^{\otimes 3}$ and
\[
N_{Aut(V)}((\mathcal{S}/\mathcal{S}_\Phi)^*)=Stab_{Aut(V)}(U\otimes U\otimes U)\cong 2^{5}(2^{10}(2^{20}:(L_5(2)\times Sym_3))).
\]
\end{prop}

\medskip

\begin{nota}\labtt{s1}
Now let $\mathcal{S}_1=\{(0,a,a)\mid a\in \Phi\}$. Then $\mathcal{S}_1$ is a totally singular subspace of $\mathcal{S}$, where $\mathcal{S}$ is defined as in \refpp{mathcals}.
Let $V'=\mathcal{V}(\mathcal{S}_1)= \oplus_{[W]\in \mathcal{S}_1} W < V$.
\end{nota}

\begin{lem}\labtt{BW}  Let $\mathcal{S}_1, V'$ be as in \refpp{s1}.
We have $ V'\cong V_{EE_8}^+ \otimes V_{BW_{16}}^+$, where $BW_{16}$ denotes the Barnes-Wall lattice of rank $16$.
\end{lem}

\pf Let $E\cong E_8$ be an overlattice of $EE_8$. Then the lattice
\[
B = span\{EE_8\perp EE_8\}\cup \{ (\a, \a)\mid \a \in E\}< E\perp E
\]
has the discriminant group $2^8$ and is isomorphic to $BW_{16}$.

Since there is only one orbit of maximal totally singular spaces of $R(U)$ under
$Aut(U)$, we may assume $$\Phi= span\{ V_{\a+EE_8}^+\mid \a+EE_8\in
E/EE_8\}\cup \{ V_{EE_8}^-\}.$$ Then we have $V'=\mathcal{V}(\mathcal{S}_1) \cong
V_{EE_8}^+\otimes V_B^+\cong V_{EE_8}^+ \otimes V_{BW_{16}}^+$ as desired.
\eop

\medskip

Recall from \cite{sh1} that $Aut(V_{BW_{16}}^+)\cong 2^{16}. \Omega^+(10,2)$.

\begin{rem}  Let $\mathcal{S}_1, V'$ be as in \refpp{s1}.
Let $\bar{\mathcal{S}}= \mathcal{S}/\mathcal{S}_1$. Then for any coset $a+ \mathcal{S}_1 \in \bar{\mathcal{S}}$, the subspace
$\mathcal{V}(a+\mathcal{S}_1)= \oplus_{[W]\in a+\mathcal{S}_1} W$ is an irreducible $V'$-module. Moreover, the dual group $\bar{\mathcal{S}}^*$ acts on $V$.
\end{rem}

\begin{prop}  \labtt{bar(s)}
 Let $\mathcal{S}_1, V'$ be as in \refpp{s1}.
Let $\bar{\mathcal{S}}= \mathcal{S}/\mathcal{S}_1$. Then
\[
N_{Aut(V)}(\bar{\mathcal{S}}^*)=Stab_{Aut(V)}(V_{EE_8}^+ \otimes V_{BW_{16}}^+)\cong 2^{10}(2^{16} \Omega^+(10,2)).
\]

\end{prop}
\pf
Let $\tilde{\mathcal{S}}= \{ \mathcal{V}(a+\mathcal{S}_1)\mid a+ \mathcal{S}_1 \in \bar{\mathcal{S}}\} \subset R(V') = R(V_{EE_8}^+) \oplus R( V_{BW_{16}}^+)$.  Then $|\tilde{\mathcal{S}}|=2^{10}$.

Recall from \cite{sh1} that both $R(V_{EE_8}^+)$ and $R( V_{BW_{16}}^+)$ are $10$-dimensional quadratic spaces of maximal Witt index.
Moreover, $Aut(V_L^+)/O_2( Aut(V_L^+))$ acts faithfully on $R(V_L^+)$ for $L=EE_8$ or $BW_{16}$ \cite[Theorem 3.20]{sh1}.

Let $\rho_1: \tilde{\mathcal{S}} \to  R(V_{BW_{16}}^+)$ and $\rho_2: \tilde{\mathcal{S}} \to  R(V_{EE_8}^+)$ be natural projections. Then both $\rho_1$ and $\rho_2$ are bijections and they determine an isomorphism $$\phi= \rho_2(\rho_1)^{-1}: R(V_{BW_{16}}^+) \to R(V_{EE_8}^+).$$

For each $g\in Aut(V_{BW_{16}}^+)\cong 2^{16}. \Omega^+(10,2)$, define
${}^{\phi}\bar{g}= \phi
\bar{g}\phi^{-1} \in \Omega^{+}(10,2)$, where $\bar{ }:2^{16}. \Omega^+(10,2)
\to \Omega^+(10,2)$ is the natural map. Then the stabilizer
$N_{\tilde{\mathcal{S}}}$ of $\tilde{\mathcal{S}}$ is given by $\{ (^\phi\bar{g}, g)\mid g\in 2^{16}. \Omega^+(10,2)\}
\cong 2^{16}. \Omega^+(10,2)$. The theorem now follows from Theorem \ref{ND*}.
\eop

\section{Centralizer of an involution}\labtt{C(z)}
In this section, we shall show that
the automorphism group of
$V=\mathcal{V}(\mathcal{S})$ \refpp{VV}
has an involution $z$ such that
$V^z\cong V_\Lambda^+$ and $C_{Aut(V)}(z)\cong 2^{1+24} Co_1$, where $\Lambda$ denotes the Leech lattice.

\medskip

\begin{nota}\labtt{S01}
Let $\Phi$, $\Psi$ and $\mathcal{S}$ be defined as in Notation \refpp{mathcals}.
Let $x \in \Phi$ be a non-zero element. We denote
\[
\begin{split}
 \mathcal{S}^0= \{ (a,b,c)\in \mathcal{S}\mid \la (a,b,c), (x,0,0)\ra =0\}, \\
 \mathcal{S}^1= \{ (a,b,c)\in \mathcal{S} \mid  \la (a,b,c), (x,0,0)\ra =1\}.
\end{split}
\]
\end{nota}

\begin{de}\labtt{z}
Let $V=\mathcal{V}(\mathcal{S})=\oplus_{[W] \in \mathcal{S}} W$ be the VOA defined
in Theorem \ref{V}. Define a linear map $z: V \to V$ by

\begin{equation}\label{defofz}
z=
\begin{cases}
1 & \text{ on } W \text{ for } [W]\in \mathcal{S}^0, \\
-1 & \text{ on } W \text{ for } [W]\in \mathcal{S}^1.
\end{cases}
\end{equation}
Then $z$ is automorphism of order 2 of $V$.
\end{de}

\begin{lem}\labtt{mcb}
Let $\mcb$ be the weight 2 subspace of $V$. Then $\dim (\mcb)=196884$.  
\end{lem}

\pf By Lemma \ref{V1}, an irreducible module in $\mathcal{S}$ has minimal weight $2$ if and only
if it has the form

1. $(a, a, 0)$, where $a\in \Phi\setminus \{0\}$ or

2. $(a+b+c, a+c, b+c)$, where $a,b\in \Phi\setminus \{0\}$, $c\in \Psi\setminus\{0\}$, $a+c,b+c$ are
non-singular
and $a+b+c$ is non-zero singular.

Note that there are $2^5-1$ non-zero vectors in $\Psi$ and $a+c$ is singular if and only if $\la
a,c\ra=1$. Thus,   there are $(2^5-1)\times 2^4\times 2^4 \times 3$ such vectors in case 2 while
there are $(2^5 -1)\times 3$ vectors in case 1. By Proposition \ref{weight}, we have
\[
 \dim ( \mcb ) = 156\times 3 + (2^5-1)\times 3 \times 8^2 + (2^5-1)\times 2^4\times 2^4 \times
3\times 8
= 196884.
\]
Note that $\dim\, U_2= 156$ and $\dim (U^3)_2= 156\times 3$.
\eop

\begin{lem}
The trace of $z$ on $\mcb$ is 276.
\end{lem}
\pf
Since $x\in\Phi\setminus \{0\}$, $(x,0,0)$ is orthogonal to all elements
of the shape $(a,a,0)$, $(0,a,a)$, $(a,a,0)$ in $\mathcal{S}$.

Moreover, $(a+b+c, a+c, b+c)$ is orthogonal to $(x,0,0)$ if and only if $\la x, c\ra =0$.
Therefore, there are $(2^4-1)\times 2^4\times 2^4 \times 3$ vectors of the form
$(a+b+c, a+c, b+c)$ in $\mathcal{S}^0$ and $2^4\times 2^4\times 2^4 \times 3$ vectors in
$\mathcal{S}^1$ that have minimal weights 2.

Thus, the trace of $z$ on $\mcb$ is
\[
 (156\times 3+ (2^5-1) \times 3\times  8^2 + (2^4-1)\times 2^4\times 2^4 \times 3\times 8)
-
 (2^4\times 2^4\times 2^4 \times 3\times 8)=276
\]
as desired.
\eop

\medskip

\begin{nota}\labtt{tV}
Let $x\in \Phi\setminus \{0\}$ be as in Notation \refpp{S01}.
Let $A$ be an irreducible $U$-module such that $[A]=x\in R(U)$ and
$M=A\otimes U\otimes U$.  Then $[M]=(x,0,0)\in R(U)^3$.
Let $V$ and $z$ be as in \refpp{z}.
Denote
\[
\tilde{V} = V^z \oplus (V^z \times_{U^3} M).
\]
Then $\tilde{V}$ is also a holomorphic framed VOA of central charge $24$ \cite{ly1}.

Note also that
\[
 V^z \times_{U^3} M = (\bigoplus_{[W] \in \mathcal{S}^0} W )\times_{U^3} M = \bigoplus_{[W]\in
\mathcal{S}^0} (W \times_{U^3} M)=
 \bigoplus_{[W]\in (x,0,0)+\mathcal{S}^0} W.
\]
Thus, we also have  $\tilde{V} = \oplus_{[W]\in \tilde{\mathcal{S}}} W$,
where $\tilde{\mathcal{S}}= \mathcal{S}^0  \cup ((x,0,0)+ \mathcal{S}^0)$.

\end{nota}

The following is easy to prove.

\begin{lem}\labtt{V124}
Let $\tilde{V}$ be defined as in Notation \refpp{tV}.

1. $\dim (\tilde{V}_1) =24$.

2. $\tilde{V}$ contains a subVOA isomorphic to $(V_{EE_8})^3$.
\end{lem}

\pf
Let $x\in \Phi\setminus \{0\}$ be as in Notation \refpp{S01}. Since $V_1=0$, we have $(V^z)_1=0$ and
thus $\tilde{V}_1 <  V^z \times_{U^3} M =\bigoplus_{[W]\in (x,0,0)+\mathcal{S}^0} W$.

By the definition of $\mathcal{S}$ \refpp{mathcals}, we know that $(0,0,0)$, $(x,x,0)$, $(x,0,x)\in
\mathcal{S}^0$. Thus, we have $(x,0,0)$, $(0,x,0)$ and $(0,0,x)$ in $(x,0,0)+
\mathcal{S}^0$ and they have minimal weight $1$.

Now let $(a+c, b+c, a+b+c), a,b\in \Phi, c\in
\Psi,$ be an element of $\mathcal{S}^0$. If $s=(x,0,0)+ (a+c, b+c, a+b+c)=(x+a+c, b+c, a+b+c)$ has
minimal weight $1$, then at least one of the coordinates must be zero; otherwise, the minimal
weight $\geq \frac{3}2$. Since $\Phi\cap \Psi=0$, we have $c=0$  and $s=(x+a,b,a+b)$ for
some $a,b\in \Phi$. Since $x+a, a$ and $a+b$ are singular and $s$ has minimal weight $1$, only one
coordinate is non-zero. Hence, $(x,0,0)$, $(0,x,0)$ and $(0,0,x)$ are the only elements in $(x,0,0)+
\mathcal{S}^0$ which have minimal weight $1$.

Moreover by
\refpp{weight},  $$\dim(A\otimes U\otimes U)_1=\dim(U\otimes A\otimes U)_1=\dim (U\otimes
U \otimes A)_1=8.$$
Hence, we have $\dim \tilde{V}_1= 8+8+8=24$.

Since $Aut(V_{EE_8}^+)$ acts transitively on non-zero singular vectors \refpp{weight}(see also
\cite{sh}),  we have
$V_{EE_8}^+\oplus A \cong V_{EE_8}^+ \oplus V_{EE_8}^-=V_{EE_8}$.

Now let $\tilde{S}'=span\{ (x,0,0), (0,x,0),(0,0,x)\} < \tilde{\mathcal{S}}$. Then
$\oplus_{[W]\in \tilde{S}'} W$ is isomorphic to $(V_{EE_8})^3$ and thus $\tilde{V}$ contains
$(V_{EE_8})^3$ as a subVOA.
\eop

\medskip

As a direct consequence, we have

\begin{coro}\labtt{VLeech}
$\tilde{V} \cong V_\Lambda$ and $V^z\cong V_\Lambda^+$ (see
\refpp{z} and  \refpp{tV}).
\end{coro}

\pf   Since $V_{EE_8^3}\cong (V_{EE_8})^3$ is a full
subVOA of $\tilde{V}$,  $\tilde{V}$ is a direct sum of irreducible
$V_{EE_8^3}$-modules. Thus by \cite{dong1}, there exists an even  lattice
$EE_8^3 < L < (EE_8^3)^* $ such that
\[
\tilde{V}= \bigoplus_{\a+ EE_8^3\in L/(EE_8^3)}  V_{\a+EE_8^3} = M(1)\otimes \CC\{L\}.
\]
Hence $\tilde{V}$ is isomorphic to the  lattice VOA $V_L$ \cite{li}. Note  that $V_{\a+EE_8^3},
\a\in L/(EE_8^3),$ are irreducible $V_{EE_8^3}$-modules and $\tilde{V}$ is a direct sum of simple
current modules of $V_{EE_8^3}$ \cite{dong1}. Therefore,
$\tilde{V}\cong V_L$ also follows from the uniqueness of simple current
extensions \cite{dm}.

Recall that
$$\dim (V_L)_1 = rank\, L + |L(2)|,$$ where $L(2)=\{\a\in L\mid \la \a,\a\ra=2\}$
\cite{flm}.

Thus $L(2)=\emptyset$ as $rank\, L=24$ and $\dim ( \tilde{V}_1) =24$.   Moreover, $L$ is unimodular
since $\tilde{V}$ is holomorphic \cite{dong1}.
Since the Leech lattice is the only even unimodular lattice with no roots, $L\cong \Lambda$
and $\tilde{V}\cong
V_\Lambda$.

Recall that $\tilde{V}= V^z \oplus (V^z \times_{U^3} M) $ \refpp{tV}. We define an  automorphism $g$ on $\tilde{V}$ by
\[
 g=
\begin{cases}
 1 & \text{ on } V^z,\\
-1 & \text{ on } V^z \times_{U^3} M.
\end{cases}
\]
Since $V_1=0$, we have $V^z_1=0$ and $\tilde{V}_1\subset V^z \times_{U^3} M$. Therefore,
$g$ acts as $-1$ on $\tilde{V}_1$ and thus it is conjugate to the automorphism $\theta$, a lift of
the $(-1)$-isometry
of $\Lambda$ by \cite{dgh}. Hence, we have
\[
 V^z= \tilde{V}^g \cong V_\Lambda^+
\]
as desired.\eop

\bigskip

Let $V^-$ be the $(-1)$-eigenspace of $z$ in $V$.  Then $V^-= \oplus_{[W]\in
\mathcal{S}^1}W$ is a direct sum of simple current modules of $U^3$. Since $z$
acts trivially on $U^3$,  $U^3\subset V^z$. Hence, by \cite[Theorem 5,4]{M2},
for any irreducible $V^z$-module $X$, we have
\[
\sum_{W} N_{V^-\  X}^{W} \leq 1
\]
where $W$ runs through all isomorphism types of irreducible $V^z$-modules.
Moreover, $V^-\times_{V^z} X\neq 0$ since $V^z$ is rational.  Thus $V^-$ is also
a simple current module of $V^z$.

Thus, by Theorem \ref{ND*}, we have

\begin{thm} \labtt{c=cz}
Let $z$ be defined as in \refpp{z}. Then we have
an exact sequence
\[
1 \to \la z\ra  \to C_{Aut(V)}(z) \overset{\eta} \to Aut(V^z)\to 1
\]
and $C_{Aut(V)}(z)$ has the shape $ 2 \cdot 2^{24} Co_1$.
\end{thm}

\pf
Let $N$ be an irreducible $V^z$-module.
Since $V$ is a framed VOA, by Theorem 1 of
\cite{ly1}, the $V^z$-module $V\times_{V^z} N$ has a structure
of an irreducible $V$-module or an irreducible $z$-twisted $V$-module.

Thus, every irreducible module
of $V^z$ can be embedded into  $V$ or the unique irreducible $z$-twisted module of $V$.
Note that $V$ is holomorphic. Therefore, it  has only one irreducible module, namely $V$ and a unique $z$-twisted module, up to isomorphisms \cite{dlm1}.

Note that $V=V^z\oplus V^-$ and $V\times_{U^3} M= (V^z\times_{U^3} M) \oplus (V^- \times_{U^3} M) $ is
the unique $z$-twisted module \cite{ly1}, where $M$ is defined as in \refpp{tV}.
Thus,
$V^z$ has exactly $4$ inequivalent irreducible modules,
namely, $V^z$, $V^-$,
$V^z\times_{U^3} M $ and $V^- \times_{U^3} M$.  Since $V_1=0$, it is clear that
the minimal weight of $V^z$ and $V^-$ are $0$ and $2$, respectively. On the
other hand, the minimal weight of $V^z\times_{U^3} M $ is $1$ and the weights
of $V^- \times_{U^3} M$ are in $\frac{1}2 +\ZZ$. Thus, $g\circ V^-\cong V^-$ for
all $g\in Aut(V^z)$ since $g\circ V^-$ and $V^-$ have the same characters. The
theorem now follows from Theorem \ref{ND*} and the fact that $Aut(V_\Lambda^+)\cong 2^{24} Co_1$
\cite{sh1}. \eop

\begin{nota}\labtt{bigc}
The group $C$ is the group which is in the middle of the short exact sequence of \refpp{c=cz}.  It has the shape
$2^{1+24} Co_1$.   (There are two 2-constrained groups of the shape $2^{1+24} Co_1$ \cite{gr76}.)
\end{nota}

\begin{rem}\labtt{vtildevnatural}
Recall from \cite{dong} that $V_\Lambda^+$ has exactly $4$ inequivalent
irreducible modules, namely, $V_\Lambda^+$, $V_\Lambda^-$,
$V_\Lambda^{T,+}$ and $V_\Lambda^{T,-}$, and their minimal weights are $0$,
$1$, $2$ and $\frac{3}2$, respectively. Since $V^z\cong V_\Lambda^+$ and
$V_1=0$, it is easy to see that $V=V^z\oplus V^-\cong V_\Lambda^+\oplus
V_\Lambda^{T,+}$ as a $V_\Lambda^+$-module.

We can furthermore say that  $V$ is isomorphic to the
Frenkel-Lepowsky-Meurman Moonshine VOA $V^\natural$, by the uniqueness of simple current
extensions \cite{dm,Hu} or by use of  \cite{ly} or  \cite{dgl} since $V$ is framed.
If we wish to avoid quoting these results, then after we prove \refpp{identification}, we can claim the more modest result  that $V$ has Moonshine type.
\end{rem}

\subsection{Conformal vectors of central charge $\frac{1}2$}\label{sec:3.1}

\begin{nota}\labtt{cvcch}
We shall use $\cvcch$ to mean simple conformal vectors of central charge $\frac{1}2$.  \end{nota}

Since $V^z\cong V_\Lambda^+$, all $\cvcch$ in the VOA $V^z$ are classified in
\cite{ls}. There are exactly two classes of $\cvcch$, $AA_1$-type and
$EE_8$-type, up to the conjugacy of $Aut(V_\Lambda^+)$. Since $\Lambda$ is
generated by norm $4$ vectors, it can be shown by \cite[Proposition 3.2]{ls} that
Miyamoto involutions associated to $AA_1$-type $\cvcch$ generate a normal
subgroup $H_\Lambda \cong Hom(\Lambda/2\Lambda, \ZZ_2) \cong 2^{24}$ in
$V_\Lambda^+$.

\begin{nota}\labtt{HL}
Let $H_\Lambda$ be the normal subgroup generated by $AA_1$-type Miyamoto
involution in $V_\Lambda^+$.
\end{nota}

\medskip

Next we shall show that $\eta^{-1}(H_\Lambda)$ is an extra-special group of
order $2^{25}$.

\begin{de}\labtt{AA1type}
For any $\a\in L(4)$, define
\[
\omega^\pm (\a) =\frac{1}{16} \a (-1)^2 \cdot  \vac  \pm
\frac{1}4( e^\a+e^{-\a}).
\]
Then $\omega^\pm (\a)$ are $\cvcch$ and we call them $\cvcch$ of $AA_1$-type.
\end{de}

\begin{lem}\label{inner}
Let $\omega^{\epsilon_1}(\a)$ and $\omega^{\epsilon_2}(\b)$ be $AA_1$-type $\cvcch$ in $V_\Lambda^+$, where $\a, \b\in \Lambda(4)$. Then,
\[
\la \omega^{\epsilon_1}(\a), \omega^{\epsilon_2}(\b)\ra =
\begin{cases}
\frac{1}{2^7} (\la \a, \b \ra)^2 & \text{ if } \a\neq \pm \b,\\
\frac{1}{4} \delta_{\epsilon_1, \epsilon_2} & \text{ if } \a=\pm \b.
\end{cases}
\]
\end{lem}

\pf  The lemma follows easily from the formulas \cite[Chapter 8]{flm} that
\[
 \la \a(-1)^2\cdot \vac, \b(-1)^2 \cdot \vac \ra= 2(\la \a, \b\ra)^2,
\]
\[
\la e^\a+e^{-\a} , e^\b +e^{-\b}\ra =
\begin{cases}
0 & \text{ if }  \a \neq \pm \b,\\
2 & \text{ if } \a =\pm \b
\end{cases}
\]
and the definition of $\omega^{\pm}(\a)$.
\eop

\begin{nota}\labtt{MI}
For a $\cvcch$, $e$, we denote by $t(e)$ the associated Miyamoto involution \cite{M95}.
\end{nota}

By Lemma \ref{inner} and Sakuma's Theorem \cite{sak}, we know that $t({\omega^{\epsilon_1}(\a)})$
commutes with $t({\omega^{\epsilon_2}(\b)})$ unless $\la \a, \b\ra=\pm 1$.

\begin{lem}\labtt{tt}
Let $\omega^{\epsilon_1}(\a)$ and $\omega^{\epsilon_2}(\b)$ be $AA_1$-type $\cvcch$ in $V_\Lambda^+$, where $\a, \b\in \Lambda(4)$. Then, as automorphisms of $V$,
\[
(\tau_{\omega^{\epsilon_1}(\a)} \tau_{\omega^{\epsilon_2}(\b)})^2 =
\begin{cases}
1 &\text{ if } \la \a, \b \ra \text{ is even},\\
z & \text{ if } \la \a, \b \ra= \pm 1.
\end{cases}
\]
\end{lem}

\pf
Since $
t({\omega^{\epsilon_1}(\a)} )t({\omega^{\epsilon_2}(\b)})
t({\omega^{\epsilon_1}(\a)}) =t\big( t({\omega^{\epsilon_1}(\a)}) \omega^{\epsilon_2}(\b)\big)$,
we have
\[
t({\omega^{\epsilon_1}(\a)})t({\omega^{\epsilon_2}(\b)})
t({\omega^{\epsilon_1}(\a)}) =
\begin{cases}
t({\omega^{-\epsilon_2}(\b)}) &\text{ if } \la \a, \b \ra=\pm 1\\
t({\omega^{\epsilon_2}(\b)}) & \text{ if } \la \a, \b \ra\text{ is even}.
\end{cases}
\]
Thus, we have $(t(\omega^{\epsilon_1}(\a)) t(\omega^{\epsilon_2}(\b)) )^2 =
t({\omega^{-\epsilon_2}(\b)})t(\omega^{\epsilon_2}(\b))=z$ by \cite[Lemma 5.14]{ls}
if $\la \a, \b \ra=\pm 1$ and $(t(\omega^{\epsilon_1}(\a)) t(\omega^{\epsilon_2}(\b)) )^2=1$ if $\la \a, \b \ra$ is even.
\eop
\medskip

As a corollary, we have
\begin{thm} \labtt{HL2}
The Miyamoto involutions $\{t({\omega^{\pm }(\a)})\mid \a \in \Lambda(4)\}$  generate a subgroup
isomorphic to $2^{1+24}$ in $Aut(V)$ and thus $C_{Aut(V)}(z)\cong 2^{1+24}.Co_1$.
\end{thm}

\pf First we note that $\Lambda$ is generated by norm 4 vectors. Moreover, $\Lambda/2\Lambda$ forms
a non-degenerate quadratic space over $\ZZ_2$ with the quadratic form $q(\a+2 \Lambda) = \frac{1}2
\la \a, \a\ra \mod 2$.

Let $T$ be the subgroup of $Aut(V)$ generates by $AA_1$-type Miyamoto
involutions. Then the restriction map $\eta$ induces an group homomorphism
$\eta: T \to H_\Lambda$ \refpp{HL} and we have an exact sequence
\[
 1 \to \la z\ra \to T \overset{\eta}\to H_\Lambda \cong Hom(\Lambda/2\Lambda, \ZZ_2) \to 1.
\]
Moreover, by Lemma \ref{tt},  we have
\[
 [t(\omega^{\epsilon_1}(\a)), t(\omega^{\epsilon_2}(\b))] = z^{\la \a, \b\ra}
\]
and thus $T=\eta^{-1}(H_\Lambda)\cong 2^{1+24}$. Note that
$z=t(\omega^+(\a))t(\omega^-(\a))\in T$. \eop

\begin{rem}\labtt{Htilde}
Let $\omega^{\epsilon}(\a)$ be an $AA_1$-type $\cvcch$ in $V_\Lambda^+< V_\Lambda$.
Then we may also consider $t(\omega^{\epsilon}(\a))$ as an automorphism of $V_\Lambda$. In this case, $t(\omega^{\epsilon}(\a))$ acts trivially on the Heisenberg part $M(1)$ and thus acts trivially on $(V_\Lambda)_1$ \cite{ls}.

Recall \refpp{HL}, \refpp{HL2}.
Now let $\tilde{H}_\Lambda< Aut(V_\Lambda)$ be the group generated by $\{
t(\omega^{\pm}(\a))\mid \a\in \Lambda(4)\}$. Then $\tilde{H}_\Lambda \cong
Hom(\Lambda, \ZZ_2) \cong 2^{24}$ and the restriction map $$ g\in
\tilde{H}_\Lambda \to g|_{V_\Lambda^+}\in H_\Lambda< Aut(V_\Lambda^+)$$
is an isomorphism from $\tilde{H}_\Lambda$ to $H_\Lambda$  \cite{ls}. By the
discussion above, $\tilde{H}_\Lambda$ acts trivially on $(V_\Lambda)_1$.
\end{rem}

\section{Analysis of the finite group $Aut(V)$}

In the first subsection, we prove that $Aut(V)$ is finite.  This involves a discussion of framed VOAs over both the real and complex field.
Then, in the second subsection, we  use finite group theory to complete our analysis of $Aut(V)$.

\subsection{Framed VOA over $\RR$}

First, we recall some facts about framed VOA over $\RR$ from \cite{M1,M}.

\begin{nota}
Let $Vir_\RR=\oplus_{n\in \ZZ} \RR L_n \oplus \RR \mathbf{c}$ be the Virasoro algebra over $\RR$.
For $c,h\in \RR$, let $L(c,h)_\RR$ be the irreducible highest weight module of $Vir_\RR$
of highest weight $h$ and central charge $c$ over $\RR$.
\end{nota}

The following results can be found in \cite{M}.

\begin{prop}[Corollary 2.2, 2.3 and Theorem 2.4 of \cite{M}]
$L(\frac{1}2,0)_\RR$ is a rational VOA, that is, all $L(\frac{1}2,0)_{\RR}$-modules   are completely reducible. Moreover, $L(\frac{1}2,0)_\RR$ has only 3 inequivalent irreducible modules, namely, $L(\frac{1}2,0)_\RR$, $L(\frac{1}2,\frac{1}2)_\RR$ and $L(\frac{1}2,\frac{1}{16})_\RR$, and
\[
L(\frac{1}2, h) \cong \CC\otimes L(\frac{1}2,h)_\RR
\]
for all $h=0, 1/2 $ or $1/16$.
\end{prop}

\begin{lem}[Lemma 2.5 of \cite{M}]\labtt{realInter}
Let $W$ be a VOA over $\RR$ and let $M^1, M^2,M^3$ be $W$-modules. Then
\[
\dim \left(I_{W} \binom{M^3}{M^1 \quad M^2} \right)\leq
\dim \left(I_{\CC\otimes W} \binom{\CC\otimes M^3}{\CC\otimes M^1 \quad \CC\otimes M^2} \right)
\].
\end{lem}

\begin{prop}[cf.  (2.5) and (2.6) of \cite{M}]
For $h_1, h_2, h_3\in \{0, 1/2, 1/16\}$, we have
\[
\dim \left(I_{L(\frac{1}2, 0)_\RR} \binom{L(\frac{1}2, h_3)_\RR}{  L(\frac{1}2, h_1)_\RR \quad L(\frac{1}2, h_2)_\RR} \right) = \dim \left(I_{L(\frac{1}2, 0)} \binom{L(\frac{1}2, h_3)}{  L(\frac{1}2, h_1) \quad L(\frac{1}2, h_2)}\right).
\]
In particular, the fusion rules for $L(\frac{1}2,0)_\RR$ over $\RR$ are exactly the same as the fusion rules for $L(\frac{1}2,0)$ over $\CC$.
\end{prop}

\begin{de}
A simple VOA $W$ over $\RR$ is framed if it contains a full subVOA $T$ isomorphic
to  $L(\frac{1}2,0)_\RR^n$.
\end{de}

\begin{nota}\labtt{MCR}
For any $\a=(\a_1, \dots, \a_n)\in \{0,1\}^n$ and  an even binary linear
code $E$, we define
\[
M^\a_\RR= \otimes_{i=1}^n L(\frac{1}2,\frac{\a_i}2)_\RR
\quad
and
\quad
M_{E,\RR}= \oplus_{\a\in E} M^\a_\RR.
\]
\end{nota}
Unlike the complex case, a simple current extension over $\RR$ is no longer unique and there is more than one VOA structure on $M_{E,\RR}$ (see for example \cite{M}). Nevertheless, the following holds.

\begin{prop}[Proposition 3.5 of \cite{M}]
Let $E$ be an even linear binary code. Then $M_{E,\RR}$ has a unique VOA structure over $\RR$
such that the invariant form on $M_{E,\RR}$ is positive definite.
\end{prop}

\begin{prop}\labtt{CandD}
Let $W$ be a framed VOA over $\RR$ such that its invariant form is positive definite.
Then there exists two binary codes $E$ and $D$ such that $D < E^\perp$ and
\[
W=\oplus_{\b\in D} W^\b,
\]
where $W^0\cong M_{E,\RR}$ and for each $\b\in D$, $W^\b$ is an irreducible $M_E$-module with the
$1/16$-word $\b$.  The $1/16$-word for an irreducible $M_{E,\RR}$-module is defined as in the
complex case \cite{M0}.
\end{prop}

\begin{coro}\labtt{fg}
Let $W$ be a framed VOA over $\RR$ such that its invariant form is positive definite. Then $W$ is
finitely generated as a VOA.
\end{coro}

\pf By Proposition, \ref{CandD}, $W$ contains a subVOA $W^0\cong M_{E, \RR}$ and $W$ is a direct
sum of finitely many irreducible $W^0$-modules.

It is clear that $M^0=\otimes_{i=1}^n L(\frac{1}2, 0)_\RR$ is generated by $n$ $\cvcch$.
Since $W^0\cong M_{E,\RR}$ is direct sum of finite many irreducible $M^0$-modules, $W^0$ is finitely
generated. Moreover, $W$ is a direct sum of finitely many $W^0$-irreducible modules. Hence $W$ is
finitely generated. \eop

\begin{nota}\labtt{framestab}
Let $W$ be a framed VOA with a positive definite invariant form over $\RR$ and $T\cong L(\frac{1}2,0)_\RR^n$ a Virasoro frame. Denote
\[
\begin{split}
Stab_{Aut(W)}(T)&=\{ g\in Aut(W)\mid g(T)=T\},\\
Stab^{pt}_{Aut(W)}(T)&=\{ g\in Aut(W)\mid g(v)=v \text{ for all } v\in T\}.
\end{split}
\]
\end{nota}

Since $L(\frac{1}2,0)$ and $L(\frac{1}2,0)_\RR$ have the same fusion rules, the following
can be proved by the same argument as in the complex case.

\begin{prop}[\cite{dgh,ly1}]\labtt{Fstab}
Let $W$ be a framed VOA with a positive definite invariant form over $\RR$. Then

1. $Stab^{pt}_{Aut(W)}(T)$ is a finite $2$-group.

2. $Stab^{pt}_{Aut(W)}(T)$ is normal in $Stab_{Aut(W)}(T)$ and
$Stab_{Aut(W)}(T)/Stab^{pt}_{Aut(W)}(T)$ is isomorphic to a subgroup of $Aut(D)$.

In particular, $Stab_{Aut(W)}(T)$ is a finite group.
\end{prop}

\begin{lem}[Theorem 5.1 of \cite{M0}] \labtt{finitecv}  
Let $W$ be a CFT type VOA over $\RR$. Suppose $W_1=0$ and the invariant form
on $W$ is positive definite. Then for any pair of distinct $\cvcch$ $e$ and $f$
in $W$, we have
\[
0\leq \la e,f\ra \leq \frac{1}{12}.
\]
In particular, $W$ has only finitely many $\cvcch$.
\end{lem}

\begin{prop}[cf. \cite{M}]\labtt{framedfinite}
Let $W$ be a framed VOA over $\RR$. Suppose the invariant form on $W$ is positive definite and $W_1=0$. Then
$Aut(W)$ is a finite group.
\end{prop}

\pf By Lemma \ref{finitecv}, $W$ has only finitely many $\cvcch$ and thus $W$
has only finitely many Virasoro frames.  By Proposition \ref{Fstab} (see also \cite{dgh}), the stabilizer of a Virasoro frame is a finite group. Hence $Aut(W)$
is finite. \eop

\begin{prop}[cf. \cite{M}]\labtt{framedfinite2}
Let $W$ be a framed VOA over $\RR$. Suppose the invariant form on $W$ is positive definite and $W_1=0$. Then
$Aut(\CC \otimes W)$ is a finite group.
\end{prop}

\pf
In this proof, $\otimes $ means $\otimes_{\RR}$.
There is a semilinear automorphism, denoted $\g$, on $\CC \otimes_{\RR} W$ which
fixes $\RR \otimes W$ and is $-1$ on $\RR \sqrt {-1}\otimes W$.

By \refpp{fg} (see also \cite{dgh}), $W$ is a finitely generated VOA. From
\cite{dgag}, $Aut(\CC \otimes W)$ is a finite dimensional algebraic group and on
it, $\g$ has finite fixed points \refpp{framedfinite}. Its corresponding action on
$Der(\CC \otimes W)$, the complex  Lie algebra  of derivations, is therefore fixed
point free and so acts as $-1$ on this complex Lie algebra.  In fact, we note that
the $(-1)$-eigenspace of $\g$ on $End( \CC \otimes W)$ may be identified with
the real subspace $\RR  \sqrt {-1} \otimes End(W)$ of $\CC \otimes
End(W)\cong End( \CC \otimes W )$.  Since this real subspace contains no
nontrivial complex subspaces, we conclude that $Der( \CC \otimes W )=0$. It
follows that the connected component of the identity $Aut(\CC \otimes W)^0$ is
0-dimensional, whence $Aut(\CC \otimes W)$ is finite.

\medskip

\begin{de}\labtt{realform}
A {\it real form} of a complex VOA $V$ is a real subspace $W$ which is closed under the VOA operations  and such that a real basis for $W$ is a complex basis for $V$.
Given a real form $W$ of $V$,
a {\it real form} of a $V$-module $M$ is a real subspace $N$ which is closed under action by $W$  and such that a real basis for $N$ is a complex basis for $M$.  We say that $N$ is a real form of $M$ with respect to the real form $W$ of $V$.
\end{de}

Next we shall show that the VOA $V$ constructed in \refpp{V} has a real form with a positive definite invariant form.

First we recall that the lattice VOA constructed in \cite{flm} can be defined over $\RR$.
Let $V_{L,\RR}= S(\hat{H}^-_\RR)\otimes \RR\{L\}$ be the real lattice VOA associated to an even positive definite lattice, where $H=\RR\otimes_\ZZ L$, $ \hat{H}^-= \oplus_{n\in \ZZ^+} H\otimes \RR t^{-n}$. As usual, we use $x(-n)$ to denote $x\otimes t^{-n}$ for $x\in H$ and $n\in \ZZ^+$.

\begin{nota}\labtt{liftof-1}
Let $\theta: V_{L,\RR} \to V_{L, \RR}$ be defined by
\[
\theta( x(-n_1)\cdots x(-n_k)\otimes e^\a) = (-1)^k x(-n_1)\cdots x(-n_k)\otimes e^{-\a}.
\]
Then $\theta$ is an automorphism of $V_{L,\RR}$, which is a lift of the $(-1)$-isometry of $L$ \cite{flm}. We shall denote the $(\pm 1)$-eigenspaces of $\theta$ on $V_{L,\RR}$ by $V_{L,\RR}^\pm$.
\end{nota}

The following result is well-known \cite{flm,M}.
\begin{prop}[cf. Proposition 2.7 of \cite{M}]\labtt{VLpd}
Let $L$ be an even positive definite lattice. Then
the real subspace
$\tilde{V}_{L,\RR}=V_{L,\RR}^+\oplus \sqrt{-1}V_{L,\RR}^-$
is a real form of $V_L$.
Furthermore,
the invariant form on $\tilde{V}_{L,\RR}$ is positive definite.
\end{prop}

\begin{nota}
Let $U_\RR= V_{EE_8,\RR}^+$ be a real form of $U$. Since all $\cvcch$ in $U$
are contained in $V_{EE_8,\RR}^+$ \cite{dmz,g156}, $U_\RR$ is a real  framed
VOA. In fact, $U_\RR \cong M_{RM(2,4), \RR}$ (see \refpp{MCR} ) since
$V_{EE_8}^+\cong M_{RM(2,4)}$, where $RM(2,4)$ is the 2nd order Reed-Muller
code of degree $4$.
\end{nota}

\begin{lem}[Lemma 3.10 of \cite{M}]\labtt{moduleRform}
Let $E$ be an even binary code.   Let $X$ be an irreducible  $M_{E,\RR}$-module.
Then $\CC\otimes X$ is an irreducible $M_E$-module.
\end{lem}

\begin{lem}
Let $M$ be an irreducible module of $U$ over $\CC$. Suppose  $[M]\in R(U)$ is a
non-zero singular element. Then $M$ has a positive definite real form.
\end{lem}

\pf Since $Aut(V_{EE_8}^+)$ acts transitively on non-zero singular vectors
\refpp{weight} (see also \cite{sh}),  there is $g\in Aut(U)$ such that $M\cong
g\circ V_{EE_8}^-$.  Recall from \cite{g156} that $Aut(U)$ is generated by
$\sigma$-involutions associated to $\cvcch$  and all $\cvcch\in U$ are contained
in $U_\RR$. Therefore,  $g$ keeps $U_\RR$ invariant and define an
automorphism on $U_\RR$, also.

By Proposition \ref{VLpd},  the invariant form on $ \tilde{V}_{EE_8,\RR} =U_\RR
\oplus (\RR\sqrt{-1} \otimes V_{EE_8,\RR}^-) $ is a positive definite.  Set $W=
\RR\sqrt{-1} \otimes V_{EE_8,\RR}^-$. Then $U_\RR\oplus (g\circ W) \cong
\tilde{V}_{EE_8,\RR} $ also has a positive definite invariant form. Moreover,
$\CC\otimes (g\circ W)\cong  g\circ V_{EE_8,\RR}^-\cong M$. Thus, $g\circ W$ is
a positive definite real form of $M$. \eop

\begin{nota} \labtt{Ws}
Let $\Phi$, $\Psi$ and $\mathcal{S}$ be as in \refpp{mathcals}.  Then
$\mathcal{V}(\Phi)=\oplus_{[M]\in \Phi} M \cong V_{E_8}$.  By Proposition
\ref{VLpd}, $\mathcal{V}(\Phi)$ has a positive definite real form $W\cong
\tilde{V}_{E_8,\RR}$.  Then $W$ is  a direct sum of irreducible $U_\RR$-modules.
Let $X$ be an irreducible  $U_\RR$-submodule of $W$.  Then,  $\CC \otimes X$ is an
irreducible $U$-module by Lemma \ref{moduleRform}. Since $\CC\otimes W=
\mathcal{V}(\Phi)$,  $[\CC\otimes X]=a $ for some $a\in \Phi$. Hence, for each
$a\in \Phi$, there exists a real $U_\RR$-module $M^a$ such that $a=
[\CC\otimes M^a]$ and $\oplus_{a\in \Phi} M^a \cong \tilde{V}_{E_8,\RR}$.
Similarly, for each $b\in \Psi$, there exists a real submodule $N^b$ such that
$b=[\CC\otimes W^b]$ and $\oplus_{a\in \Psi} N^b \cong \tilde{V}_{E_8,\RR}$.

Recall that a general element in $\mathcal{S}$ has the form $(a+b+c, a+c, b+c)$
for some $a,b\in \Phi$ and $c\in \Psi$ (cf. \refpp{mathcals} and \refpp{V1}). For
any $s=(a+b+c, a+c, b+c)\in \mathcal{S}$, we define
\[
W^s := ( M^{a+b}\times_{U_{\RR}} N^c) \otimes ( M^{a}\times_{U_{\RR}} N^c) \otimes ( M^{b}\times_{U_{\RR}} N^c)
\]
as a $(U_\RR)^3$-module. Note that $W^0\cong (U_{\RR})^3$.
\end{nota}

\begin{thm}\labtt{framedfinite3}
Let $V$ be the framed VOA constructed in \refpp{V}. Then $V$ has a real form $W$ such that the invariant form on $W$ is positive definite  and $W$ is framed. Thus, $Aut(V)$ is finite by \refpp{framedfinite2}.
\end{thm}

\pf  Let $W^s, s \in \mathcal{S}$, be defined as in \refpp{Ws}. We shall show that
$W= \oplus_{s\in \mathcal{S}} W^s$ has a real VOA structure such that the
invariant form on $W$ is positive definite.

By Lemma \refpp{realInter}, we know that all  $W^s, s\in \mathcal{S}, $ are
simple current  modules of $(U_\RR)^3$. Thus, by \cite[Theorem 5.25]{M}, it
suffices to show $\oplus_{s\in \mathcal{T}} W^s$ has a real VOA structure with a
positive definite invariant form for any 2-dimensional subspace $\mathcal{T}$ of
$\mathcal{S}$.

Let $(a+b+c, a+c, b+c)$ and $(a'+b'+c', a'+c', b'+c')$ be a basis of $\mathcal{T}$,
where $a,b,a',b'\in \Phi$ and $c,c'\in \Psi$.  Take $0\neq x\in \Phi$ such that
$x$ is orthogonal to $c, c'$. Then as in \refpp{z}, we define $z'$ on $W$ by
\[
z'=
\begin{cases}
1 & \text{ on } W^s \text{  if } \la s, (x,0,0)\ra =0, \\
-1 & \text{ on } W^s \text{  if } \la s, (x,0,0)\ra =1.
\end{cases}
\]
Then by the same argument as in \refpp{VLeech}, one can show that the fixed
point subspace $W^z$ can be embedded into $V_{\Lambda, \RR}^+$, which has
a real VOA structure with a positive definite invariant form.

Since $\mathcal{T}$ is orthogonal to $(x,0,0)$,  we have  $\oplus_{s\in
\mathcal{T}} W^s < W^z$  and $\oplus_{s\in \mathcal{T}} W^s$ has a real VOA
structure with a positive definite invariant form. \eop

\subsection{$Aut(V)$ is of Monster type}

\begin{prop}\labtt{charbyc}
Let $H$ be a finite group containing an involution $z$  so that $C_H(z)\cong C$, the group of \refpp{bigc}.  Then

(i) $H=O_{2'}(H)C$; or

(ii) $H$ is a simple group of order $\ordermonster$.
\end{prop}
\pf
\cite{tits84}.
A similar conclusion was obtained in \cite{sm}  under the additional assumption
that $z$ is conjugate in $H$ to an element of $O_2(C)\setminus \{z\}$.
That fusion assumption was verified  in the situation of \cite{grfg}.   For completeness, we give a
verification of this fusion result for $H=Aut(V)$ in an appendix \refpp{fusion}.
\eop

\begin{thm}\labtt{identification}
(i)
$Aut(V)$ is a finite simple group;

(ii) $|Aut(V)|= \ordermonster$.
\end{thm}
\pf
By \refpp{framedfinite3},
$Aut(V)$ is a finite group.
By \refpp{c=cz}, $C_{Aut(V)}(z)=C$.
By \refpp{normalizertrio}, $C$ is a proper subgroup of $Aut(V)$.

To prove (i), we quote \refpp{charbyc} or \cite{tits84}.
Observe that the structure of the group  in \refpp{normalizertrio} shows that 31 divides the order of
$Aut(V)/O_{2'}(Aut(V)$, whence the alternative (i) of \refpp{charbyc} does not apply here.

For (ii), use \refpp{charbyc}  or \cite{gms}.
\eop

\begin{rem}\labtt{aboutsimplicity} (i)
So far, determinations of the group order still depend on \cite{sm} or \cite{gms}.

(ii) Our VOA construction of the Monster has an
easy proof of finiteness  \refpp{finitecv},
whereas proof of finiteness in \cite{grfg} was more troublesome.  A short proof of finiteness, using the theory of algebraic groups, was given in \cite{tits84}.
\end{rem}

\begin{coro}\labtt{moonshinetype}
The VOA $V$ \refpp{VV}, constructed by Shimakura \cite{sh},  is of Moonshine type (see the Introduction for the definition).
\end{coro}

\appendix

\section{Appendix:  A fusion result}

The following is relevant to the alternate argument for
\refpp{charbyc} and in fact proves more about fusion.  It could be of some independent interest.

\begin{lem}\labtt{fusion}
The involution $z$ is conjugate in $Aut(V)$ to elements of $O_2(C)\setminus \la z \ra$ and to elements of $C\setminus O_2(C)$.
\end{lem}
\pf
We see this by examination of the group in \refpp{normalizertrio}.

Let $x\in \Phi$, $\mathcal{S}^0$, $\mathcal{S}^1$ and $z\in Aut(V)$ be defined
as in \refpp{S01} and \refpp{z}. Without loss, we may assume $x=[V_{EE_8}^-]$.

Recall the bilinear form on $R(U)$ from \cite{sh1} that
\begin{equation}\labtt{form}
\begin{split}
\la [V_{\frac{\a}2 +EE_8}^\pm ], [V_{\frac{\b}2 +EE_8}^\pm]\ra &= \la \a, \frac{\b}2\ra \mod 2,\\
\la [V_{EE_8}^-], [(V_{EE_8}^{T_\chi})^\pm]\ra &=1,
\end{split}
\end{equation}
where $\a,\b\in EE_8$ and $V_{EE_8}^{T_{\chi}}$ is an irreducible twisted
module $V_{EE_8}$ for some character $\chi$ of $EE_8^*/EE_8$.

Let $p_1: R(U)^3 \to R(U)$ be the natural projection to the first component.
Then for any $s\in \mathcal{S}^0$, $\la p_1(s), x\ra =0$  and by \refpp{form},
 we have $p_1(s) =[V_{\frac{\b}2 +EE_8}^\epsilon]$ for some $\b\in EE_8$.

Let $x'\in \Phi$ with $x'\neq x$. Then as in \refpp{z}, we can define an
automorphism $z'$ by
\[
z'=
\begin{cases}
1 & \text{ on $W$ if } [W]\in \mathcal{S}, \la [W], (x',0,0)\ra =0,\\
-1 & \text{ on $W$ if } [W]\in \mathcal{S}, \la [W], (x',0,0)\ra =1.
\end{cases}
\]

Again we may assume $x'=[V_{\frac{\a}2 +EE_8}^\epsilon]$ for some $\a\in
EE_8$ and $\epsilon=+$ or $-$. Then $z'$ acts on $V^z= \oplus_{[W]\in
\mathcal{S}^0} W$  and by \refpp{form}
\[
z'|_{V^z}= (-1)^{ \la\a, \frac{\b}2\ra}\quad  \text{ on $W$ with } p_1([W]) = [V_{\frac{\b}2 +EE_8}^\pm].
\]
Thus, $z'|_{V^z}\in H_\Lambda$ and $z'\in O_2(C)$, where $H_\Lambda$ is
defined as in \refpp{HL}.

By the same argument as in \refpp{V124} and \refpp{VLeech}, we also have
$V^{z'}\cong V_\Lambda^+$. Thus, by the uniqueness of simple current
extensions, there exists an automorphism $g$ that maps $V^z$ to $V^{z'}$ and
hence $z'= gzg^{-1}$.

\medskip

Next we shall show that $z$ is conjugate to an element in $C\setminus O_2(C)$.

Let $y\in \Psi$ such that $\la x,y\ra=1$. Define $z_y$ by
\[
z_y=
\begin{cases}
1 & \text{ on $W$ if } [W]\in \mathcal{S}, \la [W], (y,0,0)\ra =0,\\
-1 & \text{ on $W$ if } [W]\in \mathcal{S}, \la [W], (y,0,0)\ra =1.
\end{cases}
\]
Then we again have $V^{z_y}\cong V_\Lambda^+$ and $z_y$ is conjugate to $z$
in $Aut(V)$.

Note that $z_y$ also acts on $\tilde{V}=V^z\oplus (V^z\times_{U^3} M)\cong
V_\Lambda$ (see Notation \refpp{tV}).

Since $\la x, y\ra=1$, $z_y$ acts as $-1$
on $M$ and thus acts non-trivially on $\tilde{V}_1$. By \refpp{Htilde}, $z_y\notin
\tilde{H}_\Lambda$ and thus $z_y|_{V^z}\notin H_\Lambda$.

Therefore, $z_y\in C\setminus O_2(C)$ as desired.
 \eop

\end{document}